 \documentclass{amsart}

\usepackage{amssymb}

\newcommand{\error}{\text{error}}

\newcommand{\ZFCa}{{\operatorname{\mathsf {ZFC}}}}
\newcommand{\CH}{\operatorname{\mathsf {CH}}}
\newcommand{\reals}{{\mathbb R}}
\newcommand{\rationals}{{\mathbb Q}}
\newcommand{\rest}{{\mathord{\restriction}}}

\newcommand{\cov}{\operatorname{\mathsf  {cov}}}

\newcommand{\N}{{\mathcal N}}
\newcommand{\M}{{\mathcal M}}

\newcommand{\QED}{\hspace{0.1in} \square \vspace{0.1in}}

\newcommand{\Proof}{{\sc Proof} \hspace{0.2in}}
\newcommand{\lft}[2]{\mathopen\ifcase#1{}\oo\or
                        \big#2\or\Big#2\else\oo\fi} 
\newcommand{\rgt}[2]{\mathclose\ifcase#1{}\oo\or
                        \big#2\or\Big#2\else\oo\fi} 
\newcommand{\SM}{{\mathcal{SM}}}
\newcommand{\SN}{{\mathcal  {SN}}}

\theoremstyle{plain}
\newtheorem{theorem}{Theorem}
\theoremstyle{plain}
\newtheorem{lemma}[theorem]{Lemma}

\newtheorem{definition}[theorem]{Definition}

\begin{document}
\title{Strongly meager sets can be quite big}
\author{Tomek Bartoszynski}
\address{Department of Mathematics\\
Boise State University\\
Boise, Idaho 83725 U.S.A.}
\thanks{The first author was  partially supported by 
NSF grant DMS 9971282 and Alexander von Humboldt Foundation} 
\email{tomek@math.boisestate.edu, http://math.boisestate.edu/\char 126 tomek}
\author{Andrzej Nowik}
\curraddr{Institute of Mathematics \\
Polish Academy of Sciences \\
Abrahama 18 \\
81 -- 825 Sopot\\
Poland}
\address{University of Gda\'nsk \\
Institute of Mathematics \\
ul. Wita Stwosza 57 \\
80 -- 952 Gda\'nsk \\
Poland}
\email{matan@julia.univ.gda.pl}
\author{Tomasz Weiss}
\address{Department of Mathematics\\
WSRP, Siedlce}
\email{weiss@wsrp.siedlce.pl}

\begin{abstract}
Assume $\CH$. There exists a strongly meager set $X \subseteq
  2^\omega $ and a continuous function $F: 2^\omega \longrightarrow
  2^\omega $ such that  $F"(X)=2^\omega $.

\end{abstract}
\maketitle

\section{Introduction}
A set  $X \subseteq 2^\omega $ is strongly meager ($X \in \SM$) if for
every null set $H \in \N$, $X+H \neq \emptyset$.
A set  $X \subseteq 2^\omega $ is strongly measure zero ($X \in \SN$) if for
every meager  set $F \in \M$, $X+F \neq \emptyset$.

Let 
${\mathcal I}=\{X \subseteq 2^\omega : \forall F:2^\omega
\longrightarrow 2^\omega \text{ continuous}, F"(X) \neq 2^\omega \}.$

The following is well-known:
\begin{lemma}
\begin{enumerate}
\item $\mathcal I$ is a $\sigma$-ideal,
\item ${\mathcal I} \subseteq (s)_0$,
\item ${\mathcal I} \neq (s)_0$ (in $\ZFCa$).
\end{enumerate}
  
\end{lemma}

Notice that such a $\sigma$ -- ideal was defined and investigated 
in several papers, see for example \cite{C}.

Since strongly meager sets and strong measure zero sets are $(s)_0$ it
makes sense to ask if they are in $\mathcal I$.

It is well-known that $\SN \subseteq {\mathcal I}$.  In fact, if
$F:2^\omega \longrightarrow 2^\omega $ is a continuous function and $X
\in \SN $ then $F"(X) \in \SN$.

The purpose of this paper is to show: 
\begin{theorem}
  It is consistent with $\ZFCa$ that $\SM \not \subseteq {\mathcal I} $.
\end{theorem}
Let us remark that we also have,

\begin{theorem}[\cite{BaSh2001}]
  It is consistent with $\ZFCa$ that ${\mathcal I} = \SM$.
\end{theorem}

It is also easy to see that (under $\CH$ for example)
${\mathcal I} \not \subseteq  \SM$.

\section{Combinatorics}
The following theorem is the finitary version of the construction.
\begin{theorem}\label{7}
  For every $k \in \omega $, $\varepsilon,\delta >0 $ there exists $n
  \in \omega $ 
  such that  if $I \subseteq \omega $, $|I|>n$ then there is a
  partition $2^I=a^0 \dot{\cup} a^1$ such that  
  \begin{enumerate}
      \item if $X \subseteq 2^I$, $|X| \leq k$ then 
          $\left|\dfrac{|\bigcap_{x \in X} (a^0+x)|}{2^{|I|}}
            -\dfrac{1}{2^{|X|}}\right|<\delta $,
 \item If $U \subseteq 2^I$, and $\dfrac{|U|}{2^{|I|}}=a \geq
   \varepsilon $ then there exists a set $T_U \subseteq 2^I$,
            $\dfrac{|T_U|}{2^{|I|}} > 1-\delta $ such that 
for every $s \in T_U$
     $\left|\dfrac{| (a^0+s) \cap U|}{2^{|I|}} -
         \dfrac{a}{2}\right| 
       < \delta$.

  \end{enumerate}
\end{theorem}
\begin{proof}
(1)  This is a special case of a result proved in \cite{BaSh607}
Fix $ k, \varepsilon, \delta$,   and choose  the set $C
\subseteq 2^I$ randomly (for the moment $I$ is arbitrary). For each $s
\in 2^I$ decisions whether $s \in  
 C$ are made independently with  the probability of $s \in C$ equal to
 $1/2$.  
Thus the set $C$ is a result of a sequence of Bernoulli trials.
Note that by the Chebyshev's inequality, 
 the probability that $1/2+ \delta \geq |C|\cdot 2^{-|I|} \geq
 1/2-\delta $ approaches $1$ as $|I|$ goes to infinity.  

Let $S_n$ be the number of successes in $n$ independent Bernoulli
  trials with probability of success $p$.
We will need the following well--known fact.
\begin{theorem}[Bernstein's Inequality]\label{cor1}
  For every $ \delta  >0$, 
\begin{displaymath}P\left(\left|\frac{S_n}{n}-p\right| \geq \delta \right) 
\leq
2 e^{-n\delta^2/4}.\end{displaymath}
\end{theorem}

  Consider an arbitrary set $X \subseteq 2^I$. To simplify
  the notation denote $V = 2^I \setminus C$ and note that
$\bigcap_{s \in X} (C+s) = 2^I \setminus (V+X)$.
For a point $t \in 2^I$, $t \not\in X+V$ is equivalent to $(t+X)\cap
V=\emptyset$.
Thus  the probability that $t \not\in X+V$ is equal
to $ 2^{-|X|}$, as $t \not\in X+V$ means that $t+x \not \in V$ for $x
\in X$.

Let $G(X)$ be a
subgroup of $(2^I, +)$ generated by $X$. Since every element of $2^I$
has order $2$, it follows that $|G(X)|\leq 2^{|X|}$.

\begin{lemma}
  There are sets $\left\{U_j : j \leq {|G(X)|}\right\}$ such that: 
  \begin{enumerate}
  \item $ \forall j \ \forall s, t \in U_j \ \lft1(s \neq t
    \rightarrow s+t \not\in G(X)\rgt1)$, 
  \item $ \forall j \leq {|G(X)|} \ |U_j|=2^{|I|}/|G(X)|$,
  \item $ \forall i \neq j \ U_i \cap U_j = \emptyset$,
  \item $\bigcup_{j \leq |G(X)|} U_j = 2^I$.
  \end{enumerate}
\end{lemma}
\Proof
Choose  $U_j$'s to be disjoint selectors from the cosets
$2^I/G(X)$.~$\QED$ 

Note that if $t_1, t_2 \in U_j$ then the events $t_1 \in X+V$ and $t_2
\in X+V$ are independent since sets $t_1+X$ and $t_2+X$ are disjoint.
Consider the sets $X_j=
U_j \cap \bigcap_{s \in X} (C+s)$ for $j
\leq |G(X)|$.
The expected value of the size of this set is $ 2^{-|X|} \cdot
2^{|I|}/|G(X)|$. 
By theorem \ref{cor1} for each $j \leq |G(X)|$, 
\begin{displaymath}P\left(\left|\frac{|X_j|}{2^{|I|}/|G(X)|} - 2^{-|X|}\right|
\geq \delta \right) 
  \leq 2 
  e^{-2^{|I|-2}\delta^2/|G(X)|}.\end{displaymath}
It follows that for every $X \subseteq 2^I$ the probability that
\begin{displaymath} 2^{-|X|} - \delta \leq
  \frac{\left|\bigcap_{s \in X} (C+s)\right|}{2^{|I|}} 
\leq 2^{-|X|}+\delta \end{displaymath}
is at least 
$$1-
  2|G(X)|e^{-2^{|I|-2} \delta^2/|G(X)|} \geq 1-2^{|X|+1} 
e^{-2^{|I|-|X|-2}\delta^2}.$$
The probability that it happens for every $X$ of size $\leq k$ is at least
\begin{displaymath}1-2^{|I|\cdot (k+1)^2}\cdot 
  e^{- 2^{|I|-k-2}\delta^2} .\end{displaymath} 
If $k$ and $ \delta $ are fixed then this expression  approaches $1$ as
$|I|$ goes to infinity, 
since $\lim_{x \rightarrow \infty} P(x)e^{-x} = 0$ for any
  polynomial $P(x)$. 
It follows that for sufficiently large $|I|$ the probability that the
``random'' set $C$ has the required properties is $>0$. Thus
there exists an actual $C$ with these properties as well.

\bigskip

(2)
Let 
$$A=\{(s,t): t \in a^0+s\}$$
Note that $A$ is symmetric, that $(A)_s=\{t: (s,t) \in
A\}=(A)^s=\{t:(t,s)\in A\}$.

Let $U \subseteq 2^I$ be such that  $\dfrac{|U|}{2^{|I|}} =a\geq
\varepsilon $.
Consider sets $(A)^s=a^0+s$ for $s \in U$.  We want to know how many vertical
sections of the set $A \cap (2^I \times U)$  are of size
approximately $1/2$ relative to $U$.
Fix $s \in 2^I$ and consider  the set $\dfrac{|U \cap (A)_s|}{|U|}$. 
By the choice of $a_0$ the events $u \in a_0+s$ (equivalent to $u+s
\in a_0$) are pairwise independent with probability $1/2$. Thus, by
theorem \ref{cor1},
$$P\left(\left|\dfrac{|(A)_s \cap U|}{|U|} - \dfrac{1}{2}\right| \geq
  \delta\right)< 2e^{\frac{-{|U|}\delta^2}{4}}.$$
Probability that it holds for some $s$ is at most
$2^{|I|+1}e^{\frac{-{|U|}\delta^2}{4}}.$
Since $\dfrac{|U|}{2^{|I|}} =a\geq
\varepsilon $ we get
\begin{multline*}
P\left(\left|\dfrac{|(A)_s \cap U|}{|U|} - \dfrac{1}{2}\right| \geq
  \delta\right)=P\left(\left|\dfrac{|(A)_s \cap U|}{2^{|I|}} \cdot 
\dfrac{1}{a} - \dfrac{1}{2}\right| \geq
  \delta\right)= \\
P\left(\left|\dfrac{|(A)_s \cap U|}{2^{|I|}}
    -\dfrac{a}{2}\right|\geq a \delta \right)\leq 
2^{|I|+1}e^{\frac{-{|U|}\delta^2}{4}} \leq
2^{|I|+1}e^{\frac{-{2^{|I|}}\varepsilon \delta^2}{4}}
\end{multline*}
It follows that 
$$P\left(\left|\dfrac{|(A)_s \cap U|}{2^{|I|}}
    -\dfrac{a}{2}\right|\geq  \delta \right)\leq
2^{|I|+1}e^{\frac{-{2^{|I|}}\varepsilon \delta^2}{4}}
{\longrightarrow} 0 \text{ as } |I| \rightarrow \infty.$$ 

\end{proof}

\begin{definition}[\cite{BaSh607}]
Suppose that $I$ is a finite set. A distribution is a function $m: I
\longrightarrow \reals$ 
such that $$0 \leq m(x) \leq \frac{1}{|I|}.$$

Let $\overline{m}=\sum_{s \in I} m(s)$.
\end{definition}

To illustrate this concept suppose that $A \subseteq 2^\omega $ is a
measurable set and $n \in \omega $. Define $m$ on $2^n$ by 
$m(s)=\mu(A \cap [s])$ for $s \in 2^n$. A specific instance of this
definition that we will use often in the sequel is when $A$ is a
clopen set. In particular, if $K=I_1 \cup I_2$ and $J \subseteq 
2^K$, define distribution $m$ on $2^{I_1}$ as follows:
for $s \in 2^{I_1}$ let 
$$m(s)=\frac{|\{t \in J: s \subseteq t\}|}{2^{|K|}}.$$

The following theorem is an extension of theorem \ref{7} dealing with
distributions instead of sets.

\begin{theorem}\label{8}
  For every $k \in \omega $, $\varepsilon, \delta >0 $ there exists $n
  \in \omega $ 
  such that  if $I \subseteq \omega $, $|I|>n$ then there is a
  partition $2^I=a^0 \cup a^1$ such that  
  \begin{enumerate}
    \item if $X \subseteq 2^I$, $|X| \leq k$ then 
          $\left|\dfrac{|\bigcap_{x \in X} (a^0+x)|}{2^{|I|}}
            -\dfrac{1}{2^{|X|}}\right|<\delta $,
 \item If $m$ is a distribution on $2^I$, and $\overline{m}\geq
            \varepsilon $ then there exists a set $T_m \subseteq 2^I$,
            $\dfrac{|T_m|}{2^{|I|}} > 1-\delta $ such that 
for every $s \in T_m$,
     $\left|\sum \{m(t): t\in a^0+s \} -
         \dfrac{\overline{m}}{2}\right| 
       < \delta$.

  \end{enumerate}
\end{theorem}
\begin{proof}
This is a generalization of theorem \ref{7}. Suppose that $k,
\varepsilon, \delta $ are given and $a_0 \dot{\cup} a_1=2^I$ are as in
theorem \ref{7}.
First observe that if
$m=\dfrac{b}{2^{|I|}} \cdot \chi_U$,
where $U \subseteq 2^I$ and $\dfrac{|U|}{2^{|I|}}=a \geq
\varepsilon $, $0<b\leq 1$ and $\chi_U$ is a characteristic function
of the set $U$ then it follows immediately from
theorem \ref{7} that
for $s \in T_U$,
$$\left|\sum \{m(t): t\in a^0+s \} -
         \dfrac{\overline{m}}{2}\right| 
       < b \delta\leq \delta .$$
Next, note that if $\{U_i: i < \ell\}, \ \{b_i: i<\ell\}$ are such
that
\begin{enumerate}
\item $U_i  \cap U_j = \emptyset$ for $i \neq j$,
\item $0<b_i\leq 1$ for $i<\ell$,
\item $\dfrac{|U_i|}{2^{|I|}} \geq \varepsilon $ for every $i < \ell$,
\item $m_i=\dfrac{b_i}{2^{|I|}} \cdot \chi_{U_i}$.
\end{enumerate}
then for $m=\sum_{i<\ell} m_i$ and $s \in \bigcap_{i<\ell} T_{U_i}$ we have
$$\left|\sum \{m(t): t\in a^0+s \} -
         \dfrac{\overline{m}}{2}\right| 
       < \sum_{i<\ell}b_i \delta\leq \ell\cdot\delta \text{ \ and \
         }\dfrac{\left|\bigcap_{i<\ell} T_{U_i}\right|}{2^{|I|}} \geq
       1 -\ell \delta .  $$

Fix $k \in \omega,
\varepsilon> \delta>0 $. Apply theorem \ref{7} for $k,  \ 
\varepsilon'=\delta^2$ and
$\delta'=\delta^3$ to get $n$ and $a_0
\dot{\cup}a_1 = 2^I$.
We will show that $a_0, a_1$ satisfy the requirements of the theorem.

Consider an arbitrary distribution $m$ and let
$\ell=\dfrac{1}{\delta}$ (without loss of generality it is an integer).
Let 
$\{U_i: i < \ell\}$ be defined as
$$U_i=\left\{s \in 2^I: \dfrac{i \delta }{ 2^{|I|}}< m(s) \leq
\dfrac{(i+1) \delta }{2^{|I|}}\right\}.$$

Let $K=\left\{i: \dfrac{|U_i|}{2^{|I|}} \geq \delta^2\right\}$.
Put $m'=\sum_{i\in K}  \dfrac{i \delta}{ 2^{|I|}} \cdot \chi_{U_i}$ and
$U=\bigcup_{i<\ell} U_i$. 
Note that
for $s \in U$, and $i \in K$, $m(s)-m'(s) \leq \dfrac{\delta }{ 2^{|I|}}$ and
$|\overline{m} - \overline{m'}| \leq \delta  + \sum_{i<\ell}
\delta^2=2 \delta.$

Apply \ref{7} to each of the sets $\{U_i: i \in K\}$
to get sets $T_{U_i}$ and put $T_m = \bigcap_{i\in K}
T_{U_i}$. Clearly $\dfrac{|T_m|}{2^{|I|}} \geq 1-\ell\cdot \delta^3 =
1-\delta^2 $. 
Now, for $t \in T_m$,
$$ \sum \{m(s): s\in a^0+t   \}\leq
\sum \{m'(s): s\in a^0+t   \} + 2 \delta \leq
  \frac{\overline{m'}}{2} + \delta^3 + 2 \delta  + \delta^2 \leq
  \frac{\overline{m}}{2}+ 3\delta.$$
Lower estimate is the same and we get
$$\left|\sum \{m(s): s\in a^0+t   \} - \frac{\overline{m}}{2}\right|
\leq 3 \delta.$$
\end{proof}

\section{$\ZFCa$ result}
As a warm-up before proving the main result we will show a $\ZFCa$
result that uses only small portion of the combinatorial tools
developed above.

In order to show that $\SN \subseteq {\mathcal I}$ one could use the
following result:

\begin{theorem}[\cite{NW}]
  Suppose that $F:2^\omega \longrightarrow 2^\omega $ is a continuous
  function. There exists a set $H \in \M$ such that 
$$\forall z \in 2^\omega \ \exists y \in 2^\omega \ F^{-1}(y)
\subseteq H+z.$$
\end{theorem}

We will show that the measure analog of this theorem is false.
\begin{theorem}\label{warmup}
  There exists a continuous function $F:2^\omega \longrightarrow
  2^\omega $
such that for every set $G \in \N$,
$$\{z: \exists y \ F^{-1}(y) \subseteq G+z\}\in \N.$$
\end{theorem}
\begin{proof}
Let $\delta_n, k_n, I_n$ for $n \in \omega $ be such that  
\begin{enumerate}
\item $\delta_n = 4^{-n-3}$, $k_n=2^{n+3}$,
\item $I_n$ is chosen as  in Theorem \ref{7}(1) for $k=k_n$ and
  $\delta=\delta_n$ and
  \begin{enumerate}
  \item $|I_{n+1}| \geq 4^{n+10} 2^{|I_0 \cup \dots \cup I_{n}|}$
  \item $\min(I_{n+1}) \geq \max(I_n)$.
  \end{enumerate}
\end{enumerate}
Let $2^{I_n}=a^0_n \dot{\cup} a^1_n$ be a partition as in theorem
\ref{7}(1).

Define $F: 2^\omega \longrightarrow 2^\omega $ as 
$$F(x)(n) = i \iff x \rest I_n \in a^i_n.$$

Note that for every $x \in 2^\omega , \ F^{-1}(x)=\prod_n a_n^{x(n)}$
is a perfect set.

Suppose that $G \subseteq 2^\omega $ is null, and let $U \subseteq
2^\omega $ be an open set of measure $1/2>\varepsilon >0$ containing $G$.
We will show that 
$$\mu\left(\{z: \exists y \ F^{-1}(y) \subseteq U+z\}\right)
{\longrightarrow} 0 \text{ as }{\varepsilon \rightarrow 0} .$$

Let $$U_0^\emptyset = \left\{s \in 2^{I_0}: \dfrac{\mu([s]\cap
  U)}{\mu([s])} \geq \dfrac{3}{4}\right\},$$
and for $n>0$ and $t \in 2^{I_0\cup \dots \cup I_{n-1}}$ let
$$U_n^t = \left\{s \in 2^{I_n}: \dfrac{\mu([t^\frown s]\cap
  U)}{\mu([t^\frown s])} \geq 1-\dfrac{1}{2^{n+2}}\right\}.$$
Easy computation shows that $\mu(U^\emptyset_0) \leq 4/3 \varepsilon
\leq 2/3$.
For $t \in 2^{I_0\cup \dots \cup I_{n-1}}$ we say that $t$ is good if
$t \rest I_0 \not \in U^\emptyset_0$, and for $j \leq n-1$, 
$t \rest I_j \not\in U^{t \rest
  I_0\cup \dots \cup I_{j-1}}_j$.
In particular, if $t \in 2^{I_0\cup \dots \cup I_{n-1}}$ is good then
by induction we show that
$$\mu(U^t_n) \leq \frac{2^{n+2}}{2^{n+2}-1} \cdot \frac{\mu(U \cap
  [t])}{\mu([t])} \leq \frac{2^{n+2}}{2^{n+2}-1} \cdot
\left(1-\frac{1}{2^{n+1}}\right) \leq \frac{2^{n+2}-2}{2^{n+2}-1}.$$

For a good sequence $t \in 2^{I_0\cup \dots \cup I_{n-1}}$ let
$$Z^t_n=\{v \in 2^{I_n}: \exists i \in 2 \ a^i_n+v \subseteq
U^t_n\}.$$
By theorem \ref{7}(1), for $i=0,1$ and $X \subseteq 2^{I_n}$ of size
$n+10$, 
 $$ \frac{|X+a^i_n|}{2^{|I_n|}} = 1- \dfrac{|\bigcap_{x \in X}
   (a^i_n+x)|}{2^{|I_n|}} \geq 1 - \frac{1}{2^{n+10}} -
 \frac{1}{4^{n+3}} \geq \frac{2^{n+2}-2}{2^{n+2}-1}.$$
Therefore, $|Z^t_n| \leq n+10$.
Let $Z_n = \bigcup \{Z^t_n: t \in 2^{I_0\cup \dots \cup I_{n-1}} \text{
  is good}\}$. We get that $|Z_n| \leq (n+10)\cdot 2^{|I_0\cup \dots \cup
  I_{n-1}|}$, so in particular, $|Z_n| \cdot 2^{-|I_n|} \leq 2^{-n}$.
Let $$W_U=\{z \in 2^\omega : \exists n \ z \rest I_n \in Z_n\}.$$
Note that if $\varepsilon < 2^{-n-|I_0 \cup \dots \cup I_n|}$ then
$Z_0=Z_1=\dots = Z_n=\emptyset$. Therefore, $\mu(W_U) \longrightarrow
0$ as $ \varepsilon \longrightarrow 0$.
The following lemma finishes the proof.
\begin{lemma}
  If $F^{-1}(y) \subseteq U+z$ then $z \in W_U$.
\end{lemma}
\begin{proof}
  Suppose not. By induction
build a branch $r \in z+F^{-1}(y)$ such that  $r+z\rest I_0
\not\in (U)_0^\emptyset$, $r+z \rest I_1 \not \in (U)^{r+z \rest I_0}_1$,
etc. Since $U$ 
is open, it means that $r+z \not\in U$.
\end{proof}\end{proof}

\section{Main result}
\begin{theorem}
  Assume that $\cov(\N)=2^{\boldsymbol\aleph_0} $. There exists a set
  $X \in \SM$ and a continuous function $F:2^\omega \longrightarrow
  2^\omega $ such that $F"(X)=2^\omega $.
\end{theorem}

Suppose that we have  sequences $\{\varepsilon_n, \delta_n: n \in
\omega \}$ and a partition $\{I_n: n \in \omega\}$ such that 
\begin{enumerate}
\item $2^{|I_0\cup I_1 \cup \dots \cup I_n|} \cdot \delta_{n+1} < 
\varepsilon_{n+1}<2^{-n}$,
\item $2^{|I_0\cup I_1 \cup \dots \cup I_n|} \cdot \varepsilon_{n+1} <
  \delta_n$,
  \item $I_n $, $a^0_n \dot{\cup} a^1_n= 2^{I_n}$  is chosen for
    $\delta_n, \varepsilon_n$ as in theorem \ref{8}. 
\end{enumerate}

As in theorem \ref{warmup} define $F$ as 
$F(x)(n) = i \iff x \rest I_n \in a^i_n$.

\begin{lemma}
Suppose that a $\sigma$-ideal ${\mathcal I}$ on $2^\omega$, and 
$\sigma$-ideals ${\mathcal J}_x$ of $F^{-1}(x)$ are such that  the
following holds:
\begin{enumerate}
\item For every $G \in \N$, 
$\{z \in 2^\omega : \exists x \ (G+z) \cap F^{-1}(x) \not\in
{\mathcal 
  J}_x \} \in {\mathcal I}.$
\item for every $G \in \N$ and  $t \in 2^\omega
  $, 
$\{z: t \in (G+z)\} \in {\mathcal I}$.
\item $\forall x \ \cov({\mathcal J}_x)=\cov({\mathcal 
I})=2^{\boldsymbol\aleph_0} $.
\end{enumerate}
Then there exists a set $\in \SM$ such that $F"(X)=2^\omega $.
\end{lemma}
\begin{proof}
Let $\{G_\alpha: \alpha<c\}$ be enumeration of null sets, and
$\{t_\alpha: \alpha<c\}$ enumeration of $2^\omega $.
Build by induction sequences $\{x_\alpha, z_\alpha: \alpha < c\}$ such
that  
\begin{enumerate}
\item $x_\alpha \in F^{-1}(t_\alpha)$,
\item $\forall \beta <c \ x_\beta \not \in G_\alpha+z_\alpha $.
\end{enumerate}
Suppose that $\{x_\beta, z_\beta: \beta<\alpha\}$ are given.

Consider sets
$Z=\{z \in 2^\omega : \exists x \ (G_\alpha+z) \cap F^{-1}(x) \not\in
{\mathcal 
  J}_x \},$
and for $\beta < \alpha $,
$Z_\beta =\{z: x_\beta \in (G_\alpha+z) \}$.
Let $z_\alpha \not\in Z \cup \bigcup_{\beta<\alpha }
Z_\beta$. Next consider 
$F^{-1}(t_\alpha)$ and choose $x_\alpha \in F^{-1}(t_\alpha) \setminus
\bigcup_{\beta\leq\alpha } 
(G_\beta+z_\beta)$.
\end{proof}

In our case ${\mathcal I}=\N$ and ${\mathcal J}_x$ is the measure
ideal on $F^{-1}(x)$.
In other words, 
let $c_n= 2^{|I_n|-1}$, and let ${\mathcal J}$ be a $\sigma$-ideal of
null sets (with respect to the standard product measure) on
${\mathcal X}= \prod_{n=0}^\infty c_n$. Note that ${\mathcal X}$ is
chosen to be isomorphic 
(level by level) with $F^{-1}(x)$ for any $x$. Let ${\mathcal J}_x$ be
the copy of $\mathcal J$ on $F^{-1}(x)$. 

Specifically, define measure $\mu_x$ on $F^{-1}(x)$ as 
$\mu_x=\prod_{n \in \omega} \mu_n^{x(n)}$, where $\mu^i_n$ is a
normalized counting measure on $a^i_n$ ($i=0,1$).
Clearly, $\mu_x$ is essentially the Lebesgue measure  on $F^{-1}(x)$.

Now what we want to show is that 
\begin{lemma}\label{6}
  For every $G \in \N$,
$$\{z: \exists  x \ \mu_x(F^{-1}(x) \cap (G+z))> 0\}\in \N.$$
\end{lemma}

Before we go further let us briefly look at the nature of the
difficulties in proving the result using theorem \ref{warmup}. 
The problem is that the relation $F^{-1}(y) \not \subseteq G+z$ is not
additive. Quick analysis shows that every choice of a point $z$ that
will shift $G$ away from the set we are constructing has to fulfill
continuum many requirements.
That is why we shift from $F^{-1}(y) \not \subseteq G+z$ to $F^{-1}(y)
\cap ( G+z) \in {\mathcal J}_y$, an additive requirement.
We still have continuum many constraints, this time we need to find
$z$ such that $F^{-1}(y)
\cap ( G+z) \in {\mathcal J}_y$ for {\em every} $y$. 
 
\section{Proof of lemma \ref{6}}
Suppose that $G \subseteq 2^\omega $ is a null set.

\begin{lemma}
There are  sequences $\{K_n, K'_n:n \in \omega\}, \{J_n, J'_n: n \in
\omega\}$ such that  
\begin{enumerate}
\item $K_n$'s and $K_n'$'s are consecutive intervals that are unions
  of $I_m$'s,
\item $J_n \subseteq 2^{K_n}$, $J_n' \subseteq 2^{K_n'}$,
\item  $\dfrac{|J_n|}{2^{|K_n|}},
  \dfrac{|J_n'|}{2^{|K_n'|}} < \dfrac{1}{2^n}$,
\item $G \subseteq H_1 \cup H_2$, where
$H_1=\{x : \exists^\infty n \ x \rest K_n \in J_n\}$ and
$H_2=\{x : \exists^\infty n \ x \rest K_n' \in J_n'\}$.

\end{enumerate}  
\end{lemma}
\begin{proof}
Use the  theorem (and its proof) 2.5.7 of \cite{BJbook}.
\end{proof}
Clearly, if we show the lemma \ref{6} for $H_1$ and for $H_2$ then we
show it for $G$. Therefore, without loss of generality we can assume that 
$G=\{x : \exists^\infty n \ x \rest K_n \in J_n\}$, where 
$K_n, J_n$ are as above.
Moreover, we can assume that
$\dfrac{|J_n|}{2^{|K_n|}}=\dfrac{1}{2^n}$, since the property we are
interested in reflects downwards.

Now suppose that for $n \in \omega $,  $K_n=I_{k_n} \cup \dots \cup 
I_{k_{n+1}-1}$. 
Fix $z, x \in 2^\omega $ and let for $n \in \omega $,
$$J^{{x,z}}_n=(J_n+z\rest K_n) \cap \prod_{j=k_n}^{k_{n+1}-1}
    a^{x(j)}_j.$$
Of course $J^{x,z}_n$ depends only on $z \rest K_n$ and $x \rest [k_n,
k_{n+1})$.

It is easy to see that
$F^{-1}(x) \cap (G+z)=\{v \in F^{-1}(x): \exists^\infty n \ v\rest K_n
\in J^{x,z}_n\}$.
Since $F^{-1}(x)=\prod_n \prod_{j=k_n}^{k_{n+1}-1}
    a^{x(j)}_j$ it follows that
$$\mu_x(F^{-1}(x) \cap (G+z))=0 \iff \sum_n 
\frac{|J^{x,z}|}{\left|\prod_{j=k_n}^{k_{n+1}-1}
    a^{x(j)}_j\right|}< \infty. $$

Thus we need to find sets $T_n \subseteq 2^{K_n}$ such that
$\mu(\prod_n T_n)>0$ and such that if $z \in \rationals+\prod_n T_n$
then
$$\forall x\in 2^\omega   \ \sum_n 
\frac{|J^{x,z}|}{\left|\prod_{j=k_n}^{k_{n+1}-1}
    a^{x(j)}_j\right|}< \infty. $$

Fix $n \in \omega  $ and let $K_n=K=I_k \cup I_{k+1} \cup \cdots \cup 
I_{k+n}$, $J \subseteq
  2^K$, $\dfrac{|J|}{2^{|K|}} =2^{-n}\geq \varepsilon_k $.
It suffices to show that there exists a 
a set $T_J \subseteq 2^K$ such that
$\dfrac{|T_J|}{2^{|K|}} > 1-\delta_{k-1} $ and for every $s \in T_J$,
for every $t \in 2^{[k,k+n]}$,
$$\left|\dfrac{|(J+s) \cap \prod_{j=0}^{n}
    a^{t(j)}_{k+j}|}{|\prod_{j=0}^{n} a^{t(j)}_{k+j} |}-
  \dfrac{|J|}{2^{|K|}} \right| < \delta_{k-1} .$$

In this way for $z  \in \rationals+\prod_n T_n$, and $x \in
2^\omega $, and sufficiently large $n$, 
$$\frac{|J^{x, z}_n|}{\left|\prod_{j=k_n}^{k_{n+1}-1}
    a^{x(j)}_j\right|}=\dfrac{|(J+z) \cap \prod_{j=0}^{n}
    a^{x(j)}_{k+j}|}{|\prod_{j=0}^{n} a^{t(j)}_{k+j} |} \leq
  \dfrac{|J|}{2^{|K|}} + \delta_{k-1}.$$

It follows that to finish the construction it suffices to prove the
lemma below.

\begin{lemma}
  Suppose that $K=I_k \cup I_{k+1} \cup \cdots \cup I_{k+n}$, $J \subseteq
  2^K$, $\dfrac{|J|}{2^{|K|}} \geq \varepsilon_k $.
Then there exists a set $T_J \subseteq 2^K$ such that
$\dfrac{|T_J|}{2^{|K|}} > 1-\delta_{k-1} $ and for every $s \in T_J$,
$$\left|\dfrac{|(J+s) \cap \prod_{j=0}^{n}
    a^0_{k+j}|}{|\prod_{j=0}^{n} a^0_{k+j} |}-
  \dfrac{|J|}{2^{|K|}} \right| < \delta_{k-1} .$$
Moreover, for every $s \in T_J$, and every $t \in 2^{[k,k+n]}$,
$$\left|\dfrac{|(J+s) \cap \prod_{j=0}^{n}
    a^{t(j)}_{k+j}|}{|\prod_{j=0}^{n} a^{t(j)}_{k+j} |}-
  \dfrac{|J|}{2^{|K|}} \right| < \delta_{k-1} .$$
\end{lemma}
\begin{proof}
For $0\leq i\leq n$ define distribution $m_i$ on $2^{I_k\cup I_{k+1}
  \cup \dots \cup I_{k+i}}$ as
$$m_i(s)= \dfrac{|\{t \in J: s \subseteq t\}|}{2^{|K|}}.$$
Note that $\overline{m_i} = \dfrac{|J|}{2^{|K|}}$.
Observe that
the distribution $m_n$  coincides with $J$, that is,
$$m_n(s) = \left\{
  \begin{array}{ll}
\dfrac{1}{2^{|K|}}& \text{if } s \in J\\
0 & \text{otherwise}
  \end{array}\right. .$$

We will show by induction that for $i \leq n$, 
there exists a set $T_{m_i} \subseteq 2^{I_k \cup I_{k+1}\cup \dots
  \cup I_{k+i}}$,
$\dfrac{|T_{m_i}|}{2^{|I_k \cup I_{k+1}\cup \dots I_{k+i}|}} >
(1-\delta_k) \cdot \prod_{j< i} (1-\delta_{k+j})> 1-\delta_{k-1} $ such that 
for every $s \in T_{m_i}$,
     $$\left|\sum \left\{m_i(t): t \in \prod_{j=0}^{i} (a^0_{k+j}+s
         \rest I_{k+j}) \right\} -\dfrac{1}{2^{i + 1}}\cdot 
\dfrac{|J|}{2^{|K|}}
         \right| 
       < \frac{\delta_{k-1}}{2^i}.$$
In particular, for $i=n$, and $s \in T_{m_n}=T_J$,
$$\left|\sum \left\{m_n(t): t \in \prod_{j=0}^{n} (a^0_{k+j}+s \rest
    I_j) \right\} -\dfrac{1}{2^{{n +1}}}\cdot \dfrac{|J|}{2^{|K|}}
         \right| 
       < \frac{\delta_{k-1}}{2^{n}}.$$
The last equation means that for $s \in T_{J}$,
$$\left|\frac{\left|J \cap  \left(\prod_{j=0}^{n} a^0_{k+j} +s\rest
    I_j\right)\right| }{2^{|K|}} -\dfrac{1}{2^{n + 1}}\cdot 
\dfrac{|J|}{2^{|K|}}
         \right| 
       < \frac{\delta_{k-1}}{2^{n}}.$$
By moving $s$, and multiplying by $2^n$ we finally get,
$$\left|\frac{\left|(J+s) \cap  \prod_{j=0}^{n} a^0_{k+j}\right| 
}{\left|\prod_{j=0}^{n} a^0_{k+j}\right|} - \dfrac{|J|}{2^{|K|}}
         \right| 
       <  \delta_{k-1}.$$

Let
$$\left|\frac{\left|(J+s) \cap  \prod_{j=0}^{n} a^0_{k+j}\right| 
}{\left|\prod_{j=0}^{i} a^0_{k+j}\right|} - \dfrac{|J|}{2^{|K|}}
         \right|=\error_n . $$
We want to show that $\error_n
       <  \delta_{k-1}.$

Before we start induction note that we can shrink $J$ slightly so that
the resulting  set has the following property,
$$\forall i \leq n \ \forall s \in 2^{I_k\cup \dots \cup I_{k+i}} \
\left(m_i(s)\neq 0 \rightarrow m_i(s) \geq
  \frac{\varepsilon_{k+i}}{2^{|I_k\cup \dots \cup I_{k+i}|}}
\right).$$
By removing from $J$ all nodes (and their descendants) that do not
have this property we drop the ``measure'' by
$\varepsilon_{k+n}+\varepsilon_{k+n-1} + \dots+ \varepsilon_k < 2
\varepsilon_k < \dfrac{\delta_{k-1}}{2}$.
So let assume that $J$ has the above property and later add
$\dfrac{\delta_{k-1}}{2}$ to the error term.

The inductive proof is straightforward -- for $m_0$ we get $T_{m_0}$
immediately from theorem \ref{8}. 

Now consider 
 $m_{i+1}$. For each $t \in 2^{I_k\cup \dots \cup I_{k+i}}$ let $m^t_{i+1}$
 be the distribution on $2^{I_{k+i+1}}$ defined as
$m^t_{i+1}(s)=2^{|I_k\cup \dots \cup I_{k+i}|} m_{i+1}(t ^\frown s)$. 

Clearly, $m^t_{i+1}(s)\leq 2^{|I_k\cup \dots \cup I_{k+i}|}
\dfrac{1}{2^{|I_k\cup \dots \cup I_{k+i+1}|}}= \dfrac{1}{2^{|I_{k+i+1}|}}$
for every $s$.
Moreover, $\overline{m^t_{i+1}}=2^{|I_k\cup \dots \cup I_{k+i}|}m_i(t)$.
In particular,  shrinking $J$ as above yields, if $\overline{m^t_{i+1}}>0$
then
$\overline{m^t_{i+1}} \geq \varepsilon_{k+i}$.
For every $t \in 2^{I_k\cup \dots \cup I_{k+i}}$,
$\overline{m^t_{i+1}}>0$ apply theorem \ref{8} to get a set $T_t
\subseteq 2^{I_{k+i+1}}$ such that 
 $\dfrac{|T_t|}{2^{|I_{k+i+1}|}} \geq 1-\delta_{k+i+1}$,
and for every $s \in T_t$,
$$\left|\sum\left\{m^t_{i+1}(v): v \in a^0_{k+i+1}+s\right\} - 
\dfrac{\overline{m^t_{i+1}}}{2}\right|<\delta_{k+i+1}.$$

Let $T_{m_{i+1}}=T_{m_i} \times \bigcap \{T_t:
\overline{m^t_{i+1}}>0\}$.
Clearly,
\begin{multline*}
\dfrac{|T_{m_{i+1}}|}{2^{|I_k\cup\dots\cup I_{k+i+1}|}} = 
\dfrac{|T_{m_{i}}|}{2^{|I_k\cup\dots\cup I_{k+i}|}}\cdot
\dfrac{\left|\bigcap_t T_t\right|}{2^{|I_{k+i+1}|}} \geq\\
 \left((1-\delta_k)
\cdot \prod_{j< i} (1-\delta_{k+j})\right) \cdot \left(1-2^{|I_k\cup\dots\cup
  I_{k+i+1}|}\delta_{k+i+1}\right) \geq\\
 (1-\delta_k)
\cdot \prod_{j\leq i} (1-\delta_{k+j}) > 1-\delta_{k-1}.
\end{multline*}

Suppose that $s \in T_{m_{i+1}}$.
\begin{multline*}
\sum\left\{m_{i+1}(v): v \in \prod_{j=0}^{i+1}(a^0_{k+j}+s\rest
  I_{k+j})\right\}=\\
\sum\left\{\sum\left\{m_{i+1}(t^\frown v): v \in a^0_{k+i+1}+s \rest
    I_{k+i+1}\right\}: t \in \prod_{j=0}^{i}(a^0_{k+j}+s\rest
  I_{k+j})\right\}=\\
\sum\left\{\sum\left\{\frac{1}{2^{|I_k\cup\dots\cup I_{k+i}|}} m^t_{i+1}(v): 
v \in \prod_{j=0}^{i+1}(a^0_{k+j}+s\rest
  I_{k+j})\right\}:t \in \prod_{j=0}^{i}(a^0_{k+j}+s\rest
  I_{k+j})\right\}\leq \\
\sum\left\{\frac{1}{2^{|I_k\cup\dots\cup 
I_{k+i}|}}\dfrac{\overline{m^t_{i+1}}}{2} + \delta_{k+i+1}:  t \in 
\prod_{j=0}^{i}(a^0_{k+j}+s\rest
  I_{k+j})\right\}\leq\\
 \dfrac{1}{2}\sum\left\{\frac{1}{2^{|I_k\cup\dots\cup 
I_{k+i}|}}\overline{m^t_{i+1}} + 2\delta_{k+i+1}:  t \in 
\prod_{j=0}^{i}(a^0_{k+j}+s\rest
  I_{k+j})\right\}=\\
\dfrac{1}{2}\sum\left\{ m_i(t)+ 2\delta_{k+i+1}:  t \in 
\prod_{j=0}^{i}(a^0_{k+j}+s\rest
  I_{k+j})\right\}\leq \\
2^{|I_k\cup\dots\cup I_{k+i}|}\delta_{k+i+1} + \dfrac{1}{2}
\left(\dfrac{1}{2^i}\dfrac{|J|}{2^{|K|}} + \error_i\right) \leq 
 \dfrac{1}{2^{i+1}}\dfrac{|J|}{2^{|K|}} + 
\dfrac{\error_i}{2} + 
2^{|I_k\cup\dots\cup I_{k+i}|}\delta_{k+i+1},
\end{multline*}
where $\error_i$ is the error term given by the inductive hypothesis. 
That gives us 
$$\error_i \leq
\dfrac{\delta_k}{2^i}+\dfrac{\delta_{k+1}}{2^{i-1}} + \dots +
  \delta_{k+i} \leq \dfrac{\delta_{k-1}}{2^{i+1}},$$ so 
$$ \sum\left\{m_{i+1}(v): v \in \prod_{j=0}^{i+1}(a^0_{k+j}+s\rest
  I_{k+j})\right\} \leq \dfrac{1}{2^{i+1}}\dfrac{|J|}{2^{|K|}} +
\frac{\delta_{k-1}}{2^{i + 2}}.$$
The lower bound is similar, and we get for $s \in T_{m_{i+1}}$,
$$\left|\sum \left\{m_{i+1}(t): t \in \prod_{j=0}^{i+1} (a^0_{k+j}+s
         \rest I_{k+j}) \right\} -\dfrac{1}{2^{i+1}}\cdot \dfrac{|J|}{2^{|K|}}
         \right| 
       < \frac{\delta_{k-1}}{2^{i+2}}.$$
As before that yields the estimate 
$$\left|\frac{\left|(J+s) \cap  \prod_{j=0}^{n} a^0_{k+j}\right| 
}{\left|\prod_{j=0}^{n} a^0_{k+j}\right|} - \dfrac{|J|}{2^{|K|}}
         \right| 
       <  \frac{\delta_{k-1}}{2}.$$
Since we started by reducing the ``measure'' of $J$ by
$\dfrac{\delta_{k-1}}{2}$ we get the required estimate. 

\bigskip

Finally we will show the second part of the lemma.
We will proceed by induction on $n$.
If $n=0$ then $K=I_k$ and by the part already proved
$$\left|\dfrac{|(J+s) \cap a^0_k|}{|a^0_k|} -
\dfrac{|J|}{2^{|K|}}\right| < \delta_{k-1}.$$
Now
\begin{multline*}
|(J+s) \cap a^1_k|=|(J+s) \setminus ((J+s) \cap a^0_k)| \leq
|J|-\left(\dfrac{|J|}{2^{|K|}}- \delta_{k-1}\right) \cdot
|a^0_k|\leq\\
|J+s|-\dfrac{1}{2}|J| + \delta_{k-1} |a^0_k|=\dfrac{1}{2}|J| + \delta_{k-1} 
|a^0_k|
\end{multline*}
Thus
$$\dfrac{|(J+s) \cap a^1_k|}{2^{|K|}} \leq \frac{\frac{1}{2}|J| +
  \delta_{k-1} |a^0_k|}{2^{|K|}} =
\dfrac{1}{2}\frac{|J|}{2^{|K|}} + \frac{\delta_{k-1}}{2}.$$
 The lower estimate is similar so we have
$$\left|\dfrac{|(J+s) \cap
    a^1_k|}{2^{|K|}}-\dfrac{1}{2}\frac{|J|}{2^{|K|}}\right| <
\frac{\delta_{k-1}}{2},$$
and
$$\left|\dfrac{|(J+s) \cap
    a^1_k|}{|a^1_k}-\frac{|J|}{2^{|K|}}\right| <
\delta_{k-1}.$$
The rest of the proof is the repetition of the above argument, the
single step computed here shows that there is no difference whether we
use $a^0_j$ or $a^1_j$, the estimates do not change.
\end{proof}


\end{document}